\documentclass [12pt]{article}
\usepackage[russian]{babel}
\usepackage[cp1251]{inputenc}

\begin{document}
\begin{center}
\textbf{Uniqueness property for quasiharmonic functions}
\end{center}

\begin{center}
S.A.Imomkulov and Z. Sh. Ibragimov
\end{center}

\textbf{Abstract.} In this paper we consider class of continuous functions,
called quasiaharmonic functions, admitting best approximations by harmonic
polynomials. In this class we prove a uniqueness theorem by analogy with the
analytic functions.

\textbf{Key words.} Harmonic polynomials, quasiaharmonic functions,
polynomial approximations, N -- sets, $H - $regular compacts.

Let $K \subset {\rm R}^{n}$ be a compact set and $f(x) \in C(K)$. We denote
by

\[
l_{m} (f,K) = {\mathop {\inf }\limits_{\{q_{m} \}} }{\left\| {f(x) - q_{m}
(x)} \right\|}_{\infty }
\]

\noindent
the least deviation of the function $f$ on $K$ from harmonic polynomials of
degree $ \le m$. In paper [1] Zahariuta proved analogue of Bernshtein
theorem in the class of harmonic functions. That is, if

\begin{equation}
\label{eq1}
\overline {{\mathop {\lim }\limits_{m \to \infty } }} l_{m}^{1 / m} (f,K) <
1
\end{equation}

\noindent
then the function $f$ can be harmonically extended to some neighborhood of
the compact $K$. Conversely, if the function $f$ harmonically extends to
some neighbourhood of the compact $K$, then inequality (1) holds.

We denote by $qh(K)$ the class of functions $f(x)$ such that

\[
{\mathop {\underline {\lim } }\limits_{m \to \infty } }l_{m}^{1 / m} (f,K) <
1{\rm .}
\]

The class $qh(K)$ is called the class of quasiharmonic functions.

\textbf{Main Theorem. } Let $K \subset {\rm R}^{n}$\textit{ be a }$H - $\textit{regular compact and }$f \in qh(K)$\textit{. If the zero set }$E =
\{x \in K:\,\,\,f(x) = 0\}$\textit{ of the function f is not an }$N - $\textit{set, then }$f(x) \equiv 0$\textit{ on }$K.$

We note that in the case of quasianalytic functions theorems analogous to
our main theorem were proved in papers [2-5].

\textbf{The class of functions }$Lh_{0} (D)$\textbf{ }(see [8])\textbf{.
}Let $D$ be a domain from ${\rm R}^{n}$ and $h(D)$ be the space of harmonic
functions in $D$. We denote by $Lh_{\varepsilon } (D)$ - the minimal class
of functions, which contains all the functions of the form $\alpha \ln
{\left| {u(x)} \right|},$ $u(x) \in h(D)\,,\,\,\,\,0 < \alpha < \varepsilon
,\,\,$and closed under the operation of ``upper regulariation'', i.e., for
any family of functions $u_{\lambda } (x) \in Lh_{\varepsilon }
(D),\,\,\,\,\lambda \in \Lambda $, the function

\[
{\mathop {\overline {{\mathop {\lim }\limits_{y \to x} }} }\limits_{} }(\sup
{\left\{ {u_{\lambda } (y):\lambda \in \Lambda } \right\}})
\]

\noindent
also belongs to class $Lh_{\varepsilon } (D)$.

The union $Lh_{0} (D) = {\bigcup\limits_{\varepsilon > 0} {Lh_{\varepsilon }
(D)} }$ is called the class of $Lh_{0} - $functions.

In [1] (see also [8], [9],[10]) the author defined the following extremal
function: let $E \subset \subset D \subset {\rm R}^{n}$ be a compact set. We
fix $\varepsilon > 0$ and set

\[
\chi _{\varepsilon } (x,E,D) = {\mathop {\overline {{\mathop {\lim
}\limits_{y \to x} }} \sup {\left\{ {\alpha \ln {\left| {u(y)} \right|}:0 <
\alpha < \varepsilon ,u \in h(D),{\left\| {u} \right\|}_{E} \le 1\,,{\left\|
{u} \right\|}_{D}^{\alpha } \le e} \right\}}}\limits_{} }{\rm .}
\]

It is clear that $\chi _{\varepsilon } $ is monotonically decreasing as
$\varepsilon \to 0$ and that the following limit exists

\[
\chi _{0} (x,E,D) = {\mathop {\lim }\limits_{\varepsilon \to 0} }\chi
_{\varepsilon } (x,E,D){\rm .}
\]

Here $\chi _{0} (x,E,D)$ is called $\chi _{0} $-- measure of compact $E$
relative to the domain $D$.

As in the case of $P--$measure (see [6,7]), we have either $\chi _{0} (x,E,D)
\equiv 1$ or $\chi _{0} (x,E,D)\not { \equiv }1$ in the domain $D$. In the
first case the set $E \subset D$ for which $\chi _{0} (x,E,D) \equiv 1$ is
called the set of zero $\chi _{0} $--measure, and in the second case, the
set $E \subset D$ is called the set of nonzero $\chi _{0} $--measure.

Now we provide a lemma ``about two constant'' for the class of quasiharmonic
functions.

\textbf{Lemma 1.} (see [8],[9]). Let \textit{ }$D$\textit{ be a domain from }${\rm R}^{n}$\textit{and }$E \subset \subset D$\textit{ be a compact set. Then for any }$\alpha
\in (0,1),\,\,\varepsilon \in (0,1 - \alpha )$\textit{ and for any compact }$K \subset \subset D_{\alpha }
$\textit{ there exists a positive constant }$C = C(\alpha ,\varepsilon ,K,D)$\textit{ such that for all harmonicfunctions }$u(x)$\textit{ in }$D$\textit{ the following inequality holds:}

\begin{equation}
\label{eq2}
{\left\| {u} \right\|}_{K} \le C{\left\| {u} \right\|}_{E}^{1 - \alpha -
\varepsilon } {\left\| {u} \right\|}_{D}^{\alpha + \varepsilon } ,
\end{equation}

\textit{where }$D_{\alpha } = {\left\{ {x \in D:\,\chi _{0} (x,E,D) < \alpha }
\right\}}.$

This lemma is an analogue of the theorem ``about two constants'' for the
class of holomorphic functions (see for example [6,7]) and plays an
important role in the theory of harmonic functions.

We note that from inequality (2) it follows that if $\chi _{0} (x,E,D)\not {
\equiv }1$, then $E$ is the uniqueness set for the class of harmonic
functions in $D$.

\textbf{Definition 1 }(see [10], [11]). \textit{A compact} set\textit{ }$E \subset R^{n}$\textit{ is called }$H - $\textit{ regular at a point }$x^{0}$\textit{, if for any number }$b
> 1$\textit{ there exist numbers }$M > 0$\textit{ and }$r > 0$\textit{ such that for any harmonic polynomial }$P(x)$\textit{the following inequality holds: }

\[
{\left\| {P} \right\|}_{B(x^{0},r)} \le Mb^{\deg P}{\left\| {P(x)}
\right\|}_{E \cap B(x^{0},r)} ,
\]

\textit{where }$B(x^{0},r) = {\left\{ {x \in R^{n}:{\left| {x - x^{0}} \right|} < r}
\right\}}. $

If the compact H is regular at each of its point, then it is called $H -
$regular compact.

In [10] it is proved that if a compact $E$ is $H - $ regular at a point
$x^{0} \in E$, then for any neighbourhood $\Omega \supset \supset E$ we have

\[
\chi (x^{0},E,\Omega ) = 0{\rm .}
\]

\textbf{The }$N$\textbf{-- sets in }${\rm R}^{n}$\textbf{ }(see [8]). Let
$\vartheta _{k} (x) \in Lh_{0} (D)$ be a monotonically increasing sequence
of functions that are locally uniformly bounded from above. Consider the
limit

\[
{\mathop {\overline {\lim } }\limits_{y \to x} }{\mathop {\lim }\limits_{k
\to \infty } }\vartheta _{k} (y) = \vartheta (x){\rm ,}
\quad
x \in D{\rm .}
\]

Then everywhere in $D$ we have the inequality

\[
{\mathop {\lim }\limits_{k \to \infty } }\vartheta _{k} (x) \le \vartheta
(x){\rm .}
\]

\textbf{Definition 2.} \textbf{ }\textit{ A set }\textbf{\textit{ }}$E \subset {\rm
R}^{n}$\textbf{\textit{ }}\textit{is called}\textbf{\textit{ }}\textit{an }$N$\textit{-- set if for some open set }$D \supset E$\textit{ there exists a monotonically increasing sequence of functions }$\vartheta _{k}
(x) \in Lh_{0} (D)$\textit{ locally uniformly bounded from aboveand } \textit{such that the set }$E$\textit{ is a subset of a set of type }

\[
\{x \in D:\,\,{\mathop {\lim }\limits_{k \to \infty } }\vartheta _{k} (x) <
\vartheta (x)\},
\]

\textit{where }$\vartheta (x) = {\mathop {\overline {\lim } }\limits_{y \to x} }{\mathop
{\lim }\limits_{k \to \infty } }\vartheta _{k} (y)$, $x \in D.$

\textbf{Proposition 1 }[8]. \textit{If }$\vartheta _{k} (x) \in Lh_{0} (D)$\textit{ is sequence of functions locally uniformly bounded from above and }

\[
{\mathop {\overline {\lim } }\limits_{y \to x} }{\mathop {\overline {\lim }
}\limits_{k \to \infty } }\vartheta _{k} (y) = \vartheta (x),
\quad
x \in D,
\]

\textit{then the set }

\[
E = \{x \in D:\,\,{\mathop {\overline {\lim } }\limits_{k \to \infty }
}\vartheta _{k} (x) < \vartheta (x)\}
\]

\textit{consists of a countable union of }$N$\textit{ -- sets.}

Indeed, consider the sequence of functions

\[
w_{l,j} (x) = {\mathop {\max }\limits_{l \le k \le j} }\vartheta _{k}
(x){\rm .}
\]

Clearly, ${\mathop {\overline {\lim } }\limits_{k \to \infty } }\vartheta
_{k} (x) = {\mathop {\lim }\limits_{l \to \infty } }{\mathop {\lim
}\limits_{j \to \infty } }w_{l,j} (x)$. Since the sequence is monotonically
increasing in $j$, we have ${\mathop {\lim }\limits_{j \to \infty } }w_{l,j}
(x) \le {\mathop {\overline {\lim } }\limits_{y \to x} }{\mathop {\lim
}\limits_{j \to \infty } }w_{l,j} (y)$, $x \in D$ and the sets

\[
E_{l} = {\left\{ {x \in D:\,\,{\mathop {\lim }\limits_{j \to \infty }
}w_{l,j} (x) < {\mathop {\overline {\lim } }\limits_{y \to x} }{\mathop
{\lim }\limits_{j \to \infty } }w_{l,j} (y)} \right\}}{\rm ,}
\quad
l = 1,2,....{\rm ,}
\]

\noindent
are $N$-- sets. On the other hand, the sequenses

\[
{\mathop {\lim }\limits_{j \to \infty } }w_{l,j} (x){\rm и}
\quad
{\mathop {\overline {\lim } }\limits_{y \to x} }{\mathop {\lim }\limits_{j
\to \infty } }w_{l,j} (y){\rm ,}
\quad
l = 1,2,...,
\]

\noindent
are monotonically decreasing and

\[
{\mathop {\overline {\lim } }\limits_{k \to \infty } }\vartheta _{k} (x) =
{\mathop {\lim }\limits_{l \to \infty } }{\mathop {\lim }\limits_{j \to
\infty } }w_{l,j} (x) = \vartheta (x) = {\mathop {\lim }\limits_{l \to
\infty } }{\mathop {\overline {\lim } }\limits_{y \to x} }{\mathop {\lim
}\limits_{j \to \infty } }w_{l,j} (y){\rm ,}
\]

\[
x \in D\backslash {\bigcup\limits_{l = 1}^{\infty } {E_{l} } }{\rm .}
\]

It follows that

 $E \subset {\bigcup\limits_{l = 1}^{\infty } {E_{l} } }{\rm ,}$ i.e., $E =
{\bigcup\limits_{l = 1}^{\infty } {\left( {E_{l} \cap E} \right)} }$.

\textbf{Definition 3. }\textit{A set }$E \subset D$\textit{ is called }$Lh_{0} $\textit{ }-- \textit{polar relative to the domain }$D$\textit{ if there ecists a function }$\vartheta (x) \in Lh_{0}
(D)$\textit{ such that }$\vartheta (x)\not { \equiv } - \infty $\textit{ and }${\left. {\vartheta (x)}
\right|}_{E} = - \infty . $

We note that if $u(x) \in h(D)$, $u(x)\not { \equiv }0$ and $E \subset
{\left\{ {u(x) = 0} \right\}}$, then $E$ is $Lh_{0} $-- polar relative to
domain D.

\textbf{Proposition 2 }(see [8]). \textbf{ } \textit{Every }$Lh_{0} $--\textit{ polar set relative to domain }$D$\textit{ is contained in a countable union of }$N$--\textit{ sets.}

Indeed, let $E$ be an $Lh_{0} $-- polar set relative to domain $D$. Then by
definition there exists a function $\vartheta (x) \in Lh_{0} (D)$ such that
$\vartheta (x)\not { \equiv } - \infty $, ${\left. {\vartheta (x)}
\right|}_{E} = - \infty $. Consider a sequence of functions $\vartheta _{k}
(x) = {\textstyle{{1} \over {k}}}\vartheta (x)$. Clearly, $\vartheta _{k}
(x) \in Lh_{0} (D)$ and $\vartheta _{k} (x)\not { \equiv } - \infty $,
${\left. {\vartheta _{k} (x)} \right|}_{E} = - \infty $. Moreover, ${\mathop
{\lim }\limits_{k \to \infty } }\vartheta _{k} (x) = 0$ for almost all $x
\in D$ and ${\mathop {\lim \vartheta _{k} (x)}\limits_{k \to \infty } } = -
\infty $ for all $x \in E$. It then follows that

\[
E \subset {\left\{ {x \in D:\,\,{\mathop {\lim }\limits_{k \to \infty }
}\vartheta _{k} (x) < {\mathop {\overline {\lim } }\limits_{y \to x}
}{\mathop {\lim }\limits_{k \to \infty } }\vartheta _{k} (y)} \right\}}{\rm
.}
\]

On the othet hand, as was shown above, the set

\[
{\left\{ {x \in D:\,\,{\mathop {\lim }\limits_{k \to \infty } }\vartheta
_{k} (x) < {\mathop {\overline {\lim } }\limits_{y \to x} }{\mathop {\lim
}\limits_{k \to \infty } }\vartheta _{k} (y)} \right\}}
\]

Consists of a countable union of $N$ -- sets.

\textbf{Proof of Main Theorem. }Let $f(x) \in qh(K)$ and

\[
E = {\left\{ {x \in K:\,\,f(x) = 0} \right\}}{\rm .}
\]

By definition of the class $qh(K)$, there is a sequence of harmonic
polynomials $p_{m_{k} } (x)$ such that

\begin{equation}
\label{eq3}
{\mathop {\lim }\limits_{k \to \infty } }{\left\| {f - p_{m_{k} } }
\right\|}_{K}^{{\textstyle{{1} \over {m_{k} }}}} = d < 1{\rm .}
\end{equation}

Since ${\left. {f} \right|}_{E} = 0$, we have

\begin{equation}
\label{eq4}
{\mathop {\lim }\limits_{k \to \infty } }{\left\| {p_{m_{k} } }
\right\|}_{E}^{{\textstyle{{1} \over {m_{k} }}}} = d < 1{\rm .}
\end{equation}

Inequalities (3) and (4) imply that starting from some number $k_{0} $ for
all numbers $k \ge k_{0} $ following two inequalities hold:

\begin{equation}
\label{eq5}
{\left\| {p_{m_{k} } } \right\|}_{K} \le 1 + {\left\| {f} \right\|}{\rm ,}
\end{equation}

\begin{equation}
\label{eq6}
{\left\| {p_{m_{k} } } \right\|}_{E}^{{\textstyle{{1} \over {m_{k} }}}} < d
+ \varepsilon < 1{\rm ,}
\quad
0 < \varepsilon < 1 - d{\rm ,}
\end{equation}

Since $K$ is an $H$- regular compact, by the definition of $H$- regularity,
for any $b:\,\,1 < b < 1 / d + \varepsilon $ there are positive numbers $M$
and $\delta $ such that for a $\delta $ - heighbourhood $U_{\delta } =
{\left\{ {x:\,\,\,dist(x,K) < \delta } \right\}}$ of the compact $K$ we have
following estimate

\begin{equation}
\label{eq7}
{\left\| {P_{m_{k} } (x)} \right\|}_{U_{\delta } } \le Mb^{m_{k} }{\left\|
{P_{m_{k} } (x)} \right\|}_{K}
\end{equation}

On the other hand, since the set $E$ is not an $N$- set, we have$\chi _{0}
(x,E,U_{\delta } )\not { \equiv }1$ and using lemma 1 ``about two
constants'' we obtain that for any $\alpha \in (0,1),\,\,\beta \in (0,\,\,1
- \alpha ),\,\,\alpha + \beta < 1 / 2$, and for any open set $U:\,\,K
\subset \subset U \subset \subset U_{\delta ,\alpha } $, where $U_{\delta
,\alpha } = {\left\{ {x \in U_{\delta } :\,\chi _{0} (x,E,U_{\delta } ) <
\alpha } \right\}}$, there is a positive constant $C = C(\alpha ,\beta
,K,U_{\delta } )$ such that

\[
{\left\| {p_{m_{k} } (x)} \right\|}_{U} \le C{\left\| {p_{m_{k} } (x)}
\right\|}_{E}^{1 - \alpha - \beta } {\left\| {p_{m_{k} } (x)}
\right\|}_{D}^{\alpha + \beta }
\]

Now using estimations (5), (6) and (7) we obtain

\[
\begin{array}{l}
 {\left\| {p_{m_{k} } (x)} \right\|}_{U} \le C(d + \varepsilon )^{m_{k} (1 -
\alpha - \beta )} \cdot M^{\alpha + \beta }(1 + {\left\| {f}
\right\|})^{\alpha + \beta }b^{m_{k} (\alpha + \beta )} \le \\
 \,\,\,\,\,\,\,\,\,\,\,\,\,\,\,\,\,\,\,\,\,\,\,\,\,\,\, \le L \cdot (d +
\varepsilon )^{m_{k} (1 - \alpha - \beta )} \cdot (d + \varepsilon )^{ -
m_{k} (\alpha + \beta )} = L(d + \varepsilon )^{m_{k} (1 - 2(\alpha + \beta
))}, \\
 \end{array}
\]

\noindent
where $L = CM^{\alpha + \beta }(1 + {\left\| {f} \right\|})^{\alpha + \beta
}$.

Here $(d + \varepsilon )^{m_{k} (1 - 2(\alpha + \beta ))} \to 0,\,\,\,\,\,k
\to \infty $, since $d + \varepsilon < 1$ and $\alpha + \beta < 1 / 2$.
Therefore, ${\left\| {p_{m_{k} } (x)} \right\|}_{U} \to 0,\,\,\,\,\,\,\,k
\to \infty $, i.e., $p_{m_{k} } (x)$ converges uniformly to zero in a
neighbourhood of $U \supset \supset K$. It follows that $f(x) \equiv 0$ on
$K$. The proof is complete.

\begin{center}
\textbf{REFERENCES}
\end{center}

1. Zahariuta V.P. Inequalities for harmonic functions on spheroids
and their applications// Indiana University Mathematics Journal. --
USA, 2001, 50, №2

2. Bernstein S. N. Analytic functions of real variable, their origin
and means of generalisation. Sochineniya, Volume 1, 285-320.

3. Gonchar A.A. Kvazianaliticheskie klassi funksii, svyazannie s
nailuchshimi priblizheniyami ratsionalnimi funksiyami// Izv. A.N.
Armenia SSR 1971.- VI, № 2-3, pp. 148-159.

4. Szmuszkowiczowna H. Un theoreme sur les polynomes et son
application a la theorie des fonctions quasianalytiques// C.R.Acad.
Sci. Paris. -- 1934. Volume 198, рр. 1119- 1120.

5. Plesniak W. Quasianalytic functions of several complex
variables// Zeszyty Nauk. Uniw. Jagiell. Volume 15 (1971), рр.
135-145.

6. Sadullaev A.S. Plyurisubgarmonicheskie funksii// Itogi nauki i
tehniki. Sovremennie problemi matematiki. Fundamentalnie
napravleniya. - Moscow: VINITI, 1985.- Т. 8. - pp. 65 -- 111.

7. Sadullaev A.S. Plyurisubgarmonicheskie meri i emkosti na
kompleksnih mnogoobraziyah // UMN. - 1981. - V. 36(4). - pp. 53 -
105.

8. Imomkulov S.A., Saidov Y.R. O prodolzhenii
separatno-garmonicheskih funksiy// UzMJ. -- Tashkent, 2008, №4. pp.
89 -- 104.

9. Hecart Jean.-Marc. Ouverts d'harmonicite pour les functions
separemment harmoniques// Potential analysis. -- Netherland, 2000.
№2, рр. 115-126.

10. Hecart Jean-Marc. On Zahariuta`s extremal functions for harmonic
functions// Vietnam. J. Math. -- Springer-Verlag, 1999. Volume 27,
№1, рр. 53-59.

11. Nguyen Thanh Van., Djebbar B. Proprietes asymptotiques d'une
suite orthonormale de polynomes harmoniques// Bull. Soc. Math. 1989.
Volume 113, рр. 239-251.

\bigskip

\textit{Sevdiyor Akramovich Imomkulov}

\textit{Navoi State Pedagogical Institute (Uzbekistan)}

\textit{e-mail: sevdiyor\_i@mail.ru }

\bigskip

\textit{Zafar Shavkatovich Ibragimov}

\textit{Urgench State University (Uzbekistan)}

\textit{e-mail: z.ibragim@gmail.com }

\end{document}